%% file: main.tex
\documentclass[a4paper,12pt]{extarticle}
\usepackage{geometry}
\usepackage[utf8]{inputenc}

\providecommand{\keywords}[1]
{
  \textbf{\textit{Keywords---}} #1
}

\usepackage{subfig}

\include{commands}

\begin{document}

\title{Nonparametric estimation of a multivariate density under Kullback-Leibler loss with ISDE}
\date{}
\author{Louis Pujol\thanks{Université Paris-Saclay, CNRS, Inria, Laboratoire de Mathématiques d’Orsay, 91405, Orsay, France. louis.pujol@universite-paris-saclay.fr}\\ }
\maketitle

\begin{abstract}
  In this paper, we propose a theoretical analysis of the algorithm ISDE, introduced in previous work. From a dataset, ISDE learns a density written as a product of marginal density estimators over a partition of the features. We show that under some hypotheses, the Kullback-Leibler loss between the proper density and the output of ISDE is a bias term plus the sum of two terms which goes to zero as the number of samples goes to infinity. The rate of convergence indicates that ISDE tackles the curse of dimensionality by reducing the dimension from the one of the ambient space to the one of the biggest blocks in the partition. The constants reflect a combinatorial complexity reduction linked to the design of ISDE.
\end{abstract}


\newpage
\pagenumbering{arabic}

\section{NOTATIONS}\label{sec:notations}

Let $f$ be a density function (a nonnegative real function whose integral is equal to $1$) over $\mathbb{R}^d$. If we think of $f$ from a statistical viewpoint, it is natural to refer to the indices $\{1, \dots, d\}$ as the features.

Let $S \subset \left\{1, \dots, d\right\}$, we denote by $f_S$ the marginal density of $f$ over $S$. For all $x=(x_1, \dots, x_d) \in \mathbb{R}^d$
 \begin{equation}
    f_S(x) = \int f(x) \prod_{i \notin S} dx_i.
\end{equation}
With a slight abuse of notation, to highlight the fact that $f_S(x)$ does not depend on $(x_i)_{i \notin S}$, we write $f_S(x_S)$ instead of $f_S(x)$.

Let $k$ be an positive integer not greater than $d$. We denote by $\Sdk$ the set of all subsets of $\{1, \dots, d\}$ with cardinal not greater than $k$ and by $\Pdk$ the collection of all partitions of $\{1, \dots, d\}$ constructed with blocks in $\Sdk$. We also use the shortcuts $\Sd = \Sd^d$ and $\Pd = \Pd^d$.

\section{INTRODUCTION}\label{sec:introduction}

In a previous work (\cite{pujol2022isde}), we have introduced ISDE (Independence Structure Density Estimation). ISDE estimates a density $f$ from a set of iid realizations $X_1, \dots, X_N$ considering the Independence Structure (IS) hypothesis. This paper is devoted to a theoretical analysis of this algorithm. In this introduction, we review existing theory about IS, introduce ISDE, and set the goals of the present work.

\subsection{Curse of dimensionality and independence Structure}

\paragraph{Minimax Risk} Let $X_1, \dots, X_N$ be \emph{iid} realizations of a random variable in $\mathbb{R}^d$ admitting a density $f$. The goal of density estimation is to construct an estimator $\hat{f}$ of the density. We can measure the hardness of such an estimation task using the minimax framework. Assume that the true density belongs to some known model $\mathcal{F}$ and let $D$ be a (pseudo)distance on $\mathcal{F}$, the minimax risk is defined as follows:
\begin{equation}
    \mathcal{R}(D, \mathcal{F}) := \inf_{\hat{f}} \sup_{f \in \mathcal{F}} \esperance{D(f, \hat{f})}
\end{equation}
where the inf is taken over all measurable functions from the data to $\mathcal{F}$. More specifically, a great part of the literature on the topic deals with the asymptotic regime of $\mathcal{R}(D, \mathcal{F})$ with respect to $N$.

\paragraph{Hölder Balls} Let $\mathcal{U}$ be an open subset of $\mathbb{R}^d$ and $g : \mathcal{U} \rightarrow \mathbb{R}$ a function. Let $\gamma = (\gamma_1, \dots, \gamma_d) \in \mathbb{N}^d$ be a multiindex and let $|\gamma| = \sum_{i = 1}^d \gamma_i$ be its order. The partial differentiate operator $D^\gamma$ is defined as follows
\begin{equation}
    \ D^\gamma g = \frac{\partial^{|\gamma|} g}{\partial^{\gamma_1}_{1} \dots \partial^{\gamma_d}_{d}}.
\end{equation}

For a positive number $\beta$, if we denote by $s$ the larger integer strictly lower than $\beta$ and let $\delta = \beta - s \in (0, 1]$, $g$ belongs to the Hölder ball $\mathcal{H}(\beta, L)$ where $L$ is a positive real number if both following conditions are fulfilled

\begin{equation}
    \left\{
    \begin{array}{l}
        \underset{|\gamma | \leq s}{\max}\ \  \underset{x \in \mathcal{U}}{\sup} \  \Vert D^\gamma g(x) \Vert \leq L \\
        \underset{|\gamma | = s}{\max}\ \ \underset{x, y \in \mathcal{U}}{\sup} \  {\left| D^\gamma g(x) - D^\gamma g_i(y) \right|} \leq L{\Vert x-y\Vert ^\delta}.
    \end{array}
    \right.
\end{equation}

If $g$ is defined on a close subset $\mathcal{C}$ of $\mathbb{R}^d$, we say that $g \in \mathcal{H}(\beta, L)$ if the restriction of $g$ to the interior of $\mathcal{C}$.




\paragraph{Minimax Risk over Hölder Balls} In \cite{hasminskii1990density}, the minimax rate of this family of functions was studied considering $L_p$ distances. In particular, the result with the squared $L_2$ distance is the following
\begin{equation}
    \mathcal{R}\left( \Vert . \Vert_2^2, \mathcal{H}^\beta(d, H)\right) \sim N^{- \frac{2 \beta}{2 \beta + d}}.
\end{equation}

We can interpret this bound as a manifestation of the curse of dimensionality because of its dependence on $d$. A solution is to consider the IS model introduced in \cite{lepski2013multivariate}.

\paragraph{Independence Structure} For $k \leq d$, we define a family of functions:
\begin{equation}
    \mathcal{D}_d^k = \left\{ f: \mathbb{R}^d \rightarrow \mathbb{R}|\ \exists \Pa \in \Pdk: f(x) = \prod_{S \in \Pa} f_S(x_S) \right\}.
\end{equation}

In probabilistic terms, a density $f$ over $\mathbb{R}^d$ belongs to $\mathcal{D}_d^k$ if we can group these features into independent blocks. Another viewpoint is that the random variable characterized by $f$ admits a graphical model, a collection of disjoint fully connected cliques of size not greater than $k$. It was showed in \cite{rebelles2015mathbb} that
\begin{equation}
    \mathrm{R}\left( \Vert.  \Vert_2^2, \mathcal{H}^\beta(d, H) \cap \mathcal{D}_d^k \right) \sim N^{- \frac{2 \beta}{2 \beta + k}}.
\end{equation}

The striking fact here is that the hardness of the estimation problem is no longer related to the ambient dimension but instead to the size of the biggest block of the partition on which the density function is decomposable.

\subsection{ISDE}

As explained in \cite{pujol2022isde}, the density estimation problem under squared $L_2$ loss does not lead to a feasible algorithm. This is why we change the loss function to the Kullback-Leibler (KL) divergence. If $\hat{f}$ is an estimator of $f$, the KL loss between $f$ and $\hat{f}$ is defined as
\begin{equation}
    \KL{f}{\hat{f}} = \int \log\left( \frac{f}{\hat{f}} \right) f.
\end{equation}
This formulation is well suited to IS as it involves $\log$-densities, and a $\log$ of a product of marginal densities is a sum of $\log$ of marginal densities. From an algorithmic viewpoint, this formulation allowed us to implement an algorithm with reasonable running time and memory usage (see \cite{pujol2022isde} for details). ISDE operates as follows.

\begin{enumerate}
    \item Two independent datasets $W_1, \dots, W_m$ and $Z_1, \dots Z_n$ are extracted from $X_1, \dots, X_N$.
    \item $W_1, \dots, W_m$  is used to compute marginal density estimators $(\hat{f}_S)_{S \in \Sdk}$. Any multidimensional density estimation procedure can be used for this step.
    \item $Z_1, \dots, Z_n$ is used to compute $(\ell_n(S))_{S \in \Sdk}$ where
    \begin{equation}
        \ell_n(S) =  \frac{1}{n} \sum_{i=1}^n \log \hat{f}_S((Z_i)_S).
    \end{equation}
    \item The optimization problem
    \begin{equation}
        \max_{\Pa \in \Pdk} \sum_{S \in \Pa} \ell_n(S)
    \end{equation}
    is solved with an exact integer programming optimization procedure (branch-and-bound).
\end{enumerate}

The output is a partition $\Phat$ and an estimator of $f$ taking the form $\hat{f}_\Phat = \prod_{S \in \Phat} \hat{f}_S$. The formulation with $\log$-densities leads to a combinatorial complexity reduction. At first glance, the problem of density estimation under IS necessitates manipulating data structures of size $P_d^k$ while ISDE only lies on data structures of size $S_d^k$.

\subsection{Goal of this work}

This paper is intended to provide a theoretical analysis of ISDE by upper-bounding the quantity $\KL{f}{\hat{f}}$. In particular, we will show that the introduction of IS tackles the curse of dimensionality and that the constants in the upper-bound reflect the combinatorial complexity reduction implemented in ISDE.

\subsection{Organization of the paper}

In section \ref{sec:kldecomp} we establish a first decomposition on the risk involving oracle partitions. In section \ref{sec:conditionsobjective} we introduce the regularity conditions on the proper density and establish that an upper-bound for the uniform loss between marginal densities of $f$ and marginal estimators is sufficient to obtain a convergence result for ISDE. In section \ref{sec:unifdensity} we show that it is possible to obtain an upper-bound for uniform estimation of marginal densities for a particular estimator. Then in section \ref{sec:theorem} we state the desired upper-bound for the estimator outputted by ISDE.

\section{KULLBACK-LEIBLER RISK DECOMPOSITION}\label{sec:kldecomp}

In this section, we show that the KL loss between $f$ and $\hat{f}_\Phat$, the estimator outputted by ISDE, decomposes as the sum of three terms with a clear interpretation.

\subsection{Oracles partitions}

We denote by $\Phat$ the partition outputted by ISDE. Let $P_n[.]$ denotes the empirical measure associated with the sample $Z_1, \dots, Z_n$ and $P[.]$ the measure associated with the true density $f$. $\Phat$ is solution of the following optimisation problem :
\begin{equation}
    \Phat \in \argmin{\Pa \in \Pdk} P_n\left( - \log (\hat{f}_\mathcal{P}) \right)
\end{equation}

The partition $\Phat$ is random depending on both $W$ and $Z$. Let us define two other meaningful partitions.

\begin{align}
    \Ptilde &\in \argmin{\Pa \in \Pdk} P\left( - \log (\hat{f}_\mathcal{P}) \right) = \argmin{\Pa \in \Pdk} \KL{f}{\hat{f}_\Pa} \\
    \intertext{and}
    \Pstar &\in \argmin{\Pa \in \Pdk} P\left( - \log (f_\mathcal{P}) \right) = \argmin{\Pa \in \Pdk} \KL{f}{f_\Pa}.
\end{align}

$\Ptilde$ is a random partition depending on $W$ but not on $Z$. it is the best combination of the estimators $(\hat{f}_S)_{S \in \Sdk}$ if we consider that the quantities $\left( P(- \log \hat{f}_S) \right)_{S \in \Sdk}$ are known.

$\Pstar$ is not random. It is only a function of $k$. $f_\Pstar$ can be interpreted as the Kullback-Leibler projection of $f$ on the model $\mathcal{D}_d^k$ thanks to the following property.

\begin{myprop}{}{fstardecomp}\label{fstardecomp}
\begin{equation}
    f_\Pstar \in \argmin{g \in \mathcal{D}_d^k} \KL{f}{g}
\end{equation}
\end{myprop}

\begin{proof}
Let $g \in \mathcal{D}_d^k$ and denote by $\Pa_g$ a partition such that $g =  \prod_{S \in  {\Pa_g}} g_S$. We have
\begin{align}
    \KL{f}{g} &= \int \log\left( \frac{f}{g}\right) f \\
    &= \int \log\left( \frac{f}{f_{\Pa_g}}\right) f + \int \log\left( \frac{f_{\Pa_g}}{g}\right) f \\
    &= \KL{f}{f_{\Pa_g}} + \sum_{S \in \Pa_g} \KL{f_S}{g_S} \\
    &\leq  \KL{f}{f_{\Pa_g}}
\end{align}
with equality if $g = f_{\Pa_g}$. Then
\begin{align}
    \argmin{g \in \mathcal{D}_d^k} \KL{f}{g} &= \argmin{\Pa \in \Pdk} \min_{g \in \mathcal{D}_d^k} \KL{f}{g} \\
    &= \argmin{\Pa \in \Pdk} \KL{f}{f_{\Pa}}
\end{align}
\end{proof}

\subsection{Kullback-Leibler risk upper-bound}

We are now in a position to establish a control of the Kullback-Leibler risk for $\hat{f}_{\hat{P}}$ involving the oracles partitions.

\begin{mylem}{Kullback-Leibler risk control}{riskbound}
\begin{equation}
    \KL{f}{\hat{f}_\Phat} \leq \KL{f}{f_\Pstar} + \sum_{S_* \in \Pstar} \KL{f_{S_*}}{\hat{f}_{S_*}} + (P - P_n)(\log \hat{f}_\Ptilde - \log \hat{f}_\Phat)
\end{equation}
\end{mylem}

\begin{proof}
We start by decomposing $\KL{f}{\hat{f}_\Phat}$ as follows
\begin{align}
    \KL{f}{\hat{f}_\Phat} &= \KL{f}{f_\Pstar} \\
    &\ + \KL{f}{\hat{f}_\Pstar} - \KL{f}{f_\Pstar} \\
    &\ + \KL{f}{\hat{f}_\Ptilde} - \KL{f}{\hat{f}_\Pstar}\\
    &\ + \KL{f}{\hat{f}_\Phat} - \KL{f}{\hat{f}_\Ptilde}.
    \intertext{Then, as $\KL{f}{\hat{f}_\Ptilde} \leq \KL{f}{\hat{f}_\Pstar}$,}
    \KL{f}{\hat{f}_\Phat} &\leq \KL{f}{f_\Pstar} \ \ \ \ \\
    &\ + \KL{f}{\hat{f}_\Pstar} - \KL{f}{f_\Pstar} &&\text{(i)} \\
    &\ + \KL{f}{\hat{f}_\Phat} - \KL{f}{\hat{f}_\Ptilde} &&\text{(ii)}.
\end{align}

Now, we rewrite (i)
\begin{align}
    \KL{f}{\hat{f}_\Pstar} - \KL{f}{f_\Pstar} &= \int \log\left( \frac{f(x)}{\hat{f}_\Pstar(x)} \right) f(x) d x - \int \log\left( \frac{f(x)}{f_\Pstar(x)} \right) f(x) d x \\
    &= \int \log\left( \frac{f_\Pstar(x)}{\hat{f}_\Pstar(x)} \right) f(x) d x \\
    &= \sum_{S_* \in \Pstar} \int \log\left( \frac{f_{S_*}(x)}{\hat{f}_{S_*}(x)} \right) f_{S_*}(x) d x \\
    &= \sum_{S_* \in \Pstar} \KL{f_{S_*}}{\hat{f}_{S_*}}.
\end{align}

And we upper-bound (ii). As $P_n\left[  \log \hat{f}_\Phat - \log \hat{f}_\Ptilde \right] \geq 0$ :
\begin{align}
    \KL{f}{\hat{f}_\Phat} - \KL{f}{\hat{f}_\Ptilde} &= P\left[ -\log \hat{f}_\Ptilde \right] - P\left[ - \log \hat{f}_\Phat \right] \\
    &\leq P\left[ \log \hat{f}_\Ptilde - \log \hat{f}_\Phat \right] + P_n\left[  \log \hat{f}_\Phat - \log \hat{f}_\Ptilde \right] \\
    &= (P - P_n)(\log \hat{f}_\Ptilde - \log \hat{f}_\Phat).
\end{align}
\end{proof}

Three terms appear in the upper bound, and they can be easily interpreted.
\begin{itemize}
    \item $\KL{f}{f_\Pstar}$ is a bias term. It is the intrinsic error of the model $\mathcal{D}_d^k$ and can be thought of as a distance from $f$ to $\mathcal{D}_d^k$ thanks to proposition \ref{fstardecomp}.
    \item $\sum_{S_* \in \Pstar} \KL{f_{S_*}}{\hat{f}_{S_*}}$ is an approximation term. It is a random quantity depending on the sample $W$ and represents the error made when $f_\Pstar$ is estimated with $\hat{f}_{\Pstar}$. Conditionally to $W$, it depends on how the estimation of $\log$-likelihoods made thanks to $Z$ is accurate and quantifies our ability to output the optimal partition.
\end{itemize}

In the sequel of the paper, we will focus on upper-bounding the approximation and selection terms as they are the random quantities of interest in our problem. We treat the bias term as a structural error, and we focus on upper-bounding the quantity
\begin{equation}
    \KL{f}{\hat{f}_\Phat} - \KL{f}{f_\Pstar}.
\end{equation}
A study of the bias in a multivariate Gaussian framework can be found in appendix \ref{sec:appendixgaussian}.

\section{CONDITIONS AND OBJECTIVE}\label{sec:conditionsobjective}


\subsection{Regularity conditions}

\paragraph{Bounding condition} Density estimation under Kullback-Leibler loss is known to be a challenging problem. One work by \cite{hall1987kullback} has studied the asymptotic convergence rates for kernel estimators in a one-dimensional setting. It was shown that the tails of the kernel must be chosen appropriately regarding the tails property of the proper density to have convergent estimators. In this work, we restrict our attention to densities that are lower and upper bounded by some positive quantities. This is done to avoid hardly tractable tail behavior issues. In the sequel, we consider that the following bounding condition is valid for all $S \in \Sdk$

\begin{equation}
    \tag{BC}
       e^{-A |S|} \leq f_S \leq e^{A |S|} \ \forall S \in \Pstar \\
    \label{BC}
\end{equation}

Note that if we impose a positive lower bound on the marginal densities, we must consider that $f$ is compactly supported. In the sequel, we will suppose that the support of $f$ is $[0, 1]^d$.

\paragraph{Hölder Regularity} We will consider in the sequel that it exists $\beta \in (0, 2]$ and $L > 0$ such that $f_S \in \mathcal{H}(\beta, L)$ for all $S \in \Sdk$. We will use the following approximation property for functions in Hölder balls.

\begin{mylem}{Approximation}{controlforholder2}

Let us consider that $g \in \mathcal{H}(\beta, H)$ with $\beta \in (0, 2)$ and the domain of $g$ is $\mathcal{U} \subset \mathbb{R}^d$, then for all $x \in \mathcal{U}$ and $u$ such that $x + u \in \mathcal{U}$. If $\beta \in (0, 1]$ then
\begin{equation}
    \left| g(x) - g(x+u) \right| \leq L\Vert u \Vert^\beta.
\end{equation}
If $\beta \in (1, 2]$ then
\begin{equation}
    \left| g(x) - g(x+u) - \sum_{k=1}^d \partial_k u_k g(x)\right| \leq L\Vert u \Vert^\beta.
\end{equation}

\end{mylem}

\subsection{Objective}



Our goal is to propose an estimation procedure for the collection of marginal densities $(f_S)_{S \in \Sdk}$. If we are able to ensure, simultaneously for all $S \in \Sdk$ a uniform control 
\begin{equation}
    \tag{UC}
    \Vert \hat{f}_S - f_S \Vert_\infty \leq \epsilon_S < e^{-A |S|}(1 - e^{-A|S|}).
    \label{unifconvergence}
\end{equation}

Then we can upper-bound the approximation term and the selection term thanks to the following proposition.

\begin{myprop}{Consequences of uniform control}{csqunifconvergence}
    If the uniform control \eqref{unifconvergence} is satisfied and \eqref{BC} is true, then
    \begin{enumerate}
        \item A bounding condition is satisfied by all the estimators $(\hat{f}_S)_{S \in \Sdk}$
        \begin{equation}
            \tag{$\widehat{\text{BC}}$}
            e^{-2A|S|}\leq \hat{f}_S \leq e^{2A|S|}.
            \label{BCemp}
        \end{equation}
        \item For all $S \in \Sdk$, the Kullback-Leibler divergence between $f_S$ and $\hat{f}_S$ can be upper-bounded
        \begin{equation}
            \KL{f_S}{\hat{f}_{S}}  \leq e^{2 A |S|} \epsilon_S
        \end{equation}
        \item Conditionnaly on $W$, the selection term can be upper-bounded with high probability. More precisely if $\delta_n \in (0, 1)$ we have
        \begin{equation}
            \proba{ \left. \left| (P - P_n)(\log \hat{f}_\Ptilde - \log \hat{f}_\Phat) \right| \geq 2d\sqrt{\frac{2Ak}{n}}\sqrt{\log\left(\frac{2 S_d^k}{\delta}\right)} \right| W} \leq \delta_n.
        \end{equation}
    \end{enumerate}
\end{myprop}

\begin{proof}

\textit{Proof of 1}: Under \eqref{BC} we have
\begin{equation}
    e^{-A|S|} - \Vert \hat{f}_S - f_S \Vert_\infty \leq \hat{f}_S \leq e^{A|S|} + \Vert \hat{f}_S - f_S \Vert_\infty.
\end{equation}
Now, under \eqref{unifconvergence}
\begin{equation}
    e^{-A|S|} - \Vert \hat{f}_S - f_S \Vert_\infty \geq e^{-A|S|} - e^{-A|S|}(1 - e^{-A|S|}) = e^{-2A|S|}
\end{equation}
and
\begin{align}
    e^{A|S|} + \Vert \hat{f}_S - f_S \Vert_\infty &\leq e^{A|S|} + e^{-A|S|}(1 - e^{-A|S|}) \\
    &\leq e^{A|S|} + e^{3A|S|}\left( e^{-A|S|}(1 - e^{-A|S|}) \right) \\
    &= e^{A|S|} + e^{2A|S|}(1 - e^{-A|S|}) = e^{2A|S|}.
\end{align}

\textit{Proof of 2}: Let us compute
\begin{align}
    \KL{f_S}{\hat{f}_{S}} &= \int \log\left( \frac{f_S}{\hat{f}_{S}} \right) f_S \\
    &\leq \int \left( \frac{f_S -\hat{f}_{S}}{\hat{f}_{S}} \right) f_S \\
    \intertext{Using \eqref{BCemp}, 1 / $\hat{f}_S \leq e^{2A|S|}$}
    &\leq e^{2A|S|} \Vert f_S -\hat{f}_{S, h_m} \Vert_\infty \\
    &\leq e^{2A|S|} \epsilon_S
\end{align}

\textit{Proof of 3}:  Let $S \in \Sdk$, under \eqref{BCemp} we have $\log  \hat{f}_S \in [-2A|S|, 2A|S|]$. Using Hoeffding inequality, we obtain
\begin{equation}
    \proba{ \left| (P - P_n) \log \hat{f}_S \right| \geq \sqrt{\frac{2A|S|}{n}} \sqrt{\log \frac{2 S_d^k}{\delta}} | W} \leq \frac{\delta}{S_d^k}.
\end{equation}
Now, by union bound :
\begin{equation}
    \proba{ \sup_{S \in \Sdk}\left| (P - P_n) \log \hat{f}_S \right| \geq \sqrt{\frac{2A|S|}{n}} \sqrt{\log \frac{2 S_d^k}{\delta}} | W} \leq \delta.
\end{equation}

This leads to :
\begin{equation}
    \proba{ 2d \sup_{S \in \Sdk}\left| (P - P_n) \log \hat{f}_S \right| \leq 2d \sqrt{\frac{2Ak}{n}} \sqrt{\log \frac{2 S_d^k}{\delta}} | W} \geq 1 - \delta.
\end{equation}

Now, we remark that
\begin{align}
    \label{decompsumselectionterm}
    |(P-P_n) \log \hat{f}_{\Ptilde} - \log \hat{f}_{\Phat}| &= \left| \sum_{S \in \Ptilde} (P - P_n) \log \hat{f}_S - \sum_{S \in \Phat} (P - P_n) \log \hat{f}_S \right| \\
    &\leq \sum_{S \in \Ptilde} \left| (P-P_n) \log \hat{f}_S \right| + \sum_{S \in \Phat} \left| (P-P_n) \log \hat{f}_S \right| \\
    &\leq 2d \sup_{S \in \Sdk}\left| (P - P_n) \log \hat{f}_S \right|.
\end{align}

Then, we have
\begin{equation}
    \proba{ |(P-P_n) \log \hat{f}_{\Ptilde} - \log \hat{f}_{\Phat}| \geq 2d \sqrt{\frac{2Ak}{n}} \sqrt{\log \frac{2 S_d^k}{\delta}} | W} \leq \delta.
\end{equation}



\end{proof}

\section{UNIFORM DENSITY ESTIMATION \\ FOR MARGINAL DENSITIES}\label{sec:unifdensity}

\subsection{For a fixed $S$}

In this subsection, we fix a subset of variables $S \in \Sdk$ and we study the problem of constructing an estimator $\hat{f}_S$ based on the sample $W_1, \dots, W_m$ giving a control of $\Vert f_S - \hat{f}_S \Vert_\infty$ in order to verify \eqref{unifconvergence}. We decompose the error as a sum of a bias and a variance term as follows

\begin{equation}
\Vert f_S - \hat{f}_S \Vert_\infty \leq \underbrace{\left\Vert f_S - \esperance{\hat{f}_S} \right\Vert_\infty}_{\text{Bias}} + \underbrace{\left\Vert \esperance{\hat{f}_S} - \hat{f}_S \right\Vert_\infty}_{\text{Variance}}.
\end{equation}

\subsubsection{Bias upper-bound}

\paragraph{Choice of the kernel function} In the following we will use density estimator based on an ancillary function $K$ called kernel. $K$ is a nonnegative integrable function on $\mathbb{R}$ such that $\int K(x)dx = 1$, we consider the following assumptions on $K$:
\begin{equation}
    \tag{A.K}
    \left\{
    \begin{array}{l}
        \forall x \in \mathbb{R}, K(-x) = K(x) \\
        \text{Supp}(K) \in [-1, 1] \\
        \int x K(x) dx = 0 \\
        \Vert K \Vert_\infty < \infty
    \end{array}
    \right.
    \label{assumptionKernel}
\end{equation}

We will also assume that, if $K^{S}_{h, x} : u \mapsto \frac{1}{h^{|S|}} \prod_{k \in S} K\left( \frac{x_k - u_k}{h} \right)$, the family of function
\begin{equation}
    \mathcal{F}_S = \left\{ K^{S}_{h, x}, h > 0, x \in \mathbb{R}^{|S|} \right\}
\end{equation}
is a bounded VC class of functions. It means that it exists positive numbers $A$ and $\nu$ such that for any probability measure $P$ over $\mathbb{R}^{|S|}$ and any $\tau \in (0, 1)$ we have
\begin{equation}
    \mathcal{N}\left( \mathcal{F}_S, L_2(P) , \tau \right) \leq \left( \frac{A \Vert K \Vert_\infty}{\tau} \right)^\nu
\end{equation}
where $\mathcal{N}\left( \mathcal{F}_S, L_2(P) , \tau \right)$ is the $\tau$-covering number of $\mathcal{F}_S$ for the $L_2(P)$ distance. As proved in \cite{gine2002rates} this condition is met for almost all classical kernels. An example of kernel function $K$ satisfying all the assumptions is the Epanechnikov kernel $K_\text{Epa}$
\begin{equation}
    K_\text{Epa}(x) =  \frac{3}{4}(1-x^2) \mathbb{1}_{[-1, 1]}(x).
\end{equation}

\paragraph{Boundary issue} One must be aware of the issue induced by the fact that $f_S$ is supported on $[0, 1]^{|S|}$. We define the usual kernel density estimator (KDE) as follows. Let $h$ be a positive real number. The KDE for the marginal density $f_S$ associated with the kernel $K$, the bandwidth $h > 0$,  and the sample $W$ is defined as 

\begin{equation}
    \hat{f}_{h ,S}^{\text{KDE}}(x) = \frac{1}{mh^{|S|}} \sum_{i=1}^m  \prod_{k \in S} K\left( \frac{(W_i)_k - x_k}{h} \right)
\end{equation}

We remark that even in the samples $W_1, \dots, W_m$ belong to $[0, 1]^d$, there is no reason to have $\hat{f}_{h ,S}^{\text{KDE}}$ supported in $[0, 1]^{|S|}$.

In this setting, the bias of the classical KDE does not go to zero as $h \rightarrow 0$, illustrating the boundary issue induced by estimating a compactly supported density.

\begin{myprop}{Boundary issue}{boundaryissueforKDE}
    Let $f_S \in \mathcal{H}(2, L)$, the bias
    
    \begin{equation}
        \left\Vert \esperance{\hat{f}^{\text{KDE}}_{h, S}} - f_S \right\Vert_\infty
    \end{equation}
    
    does not tend to $0$ as $h \rightarrow 0$. 
\end{myprop}

\begin{proof}
\begin{align}
    \esperance{\hat{f}^{\text{KDE}}_{h, S}}(0) - f_S(0) &= \frac{1}{h^{|S|}} \int_{[0, 1]^{|S|}}  f_S(t) \prod_{k \in S}K\left( \frac{t_k}{h}\right) dt - f_S(0)\\
    &= \int_{[-1, 1]^{|S|}}  [f_S(hu) - f_S(0)] \prod_{k \in S} K(u_k) du_k \\
    \intertext{As $\int xK(x) = 0$}
    &= \int_{[-1, 1]^{|S|}} \left[f_S(hu) - f_S(0) - h \sum_{k \in S} u_k \partial_k f_S(0) \right] \prod_{k \in S} K(u_k) du_k \\
    \intertext{Now as $f_S(x) = 0$ for $x \notin [0, 1]^{|S|}$}
    &= \int_{[0, 1]^{|S|}} \left[f_S(hu) - f_S(0) - h \sum_{k \in S} u_k \partial_k f_S(0) \right] \prod_{k \in S} K(u_k) du_k \notag \\
    &\ \ - f_S(0) \int_{[-1, 1]^{|S|}\setminus[0, 1]^{|S|}} \prod_{k \in S} K(u_k) du_k \notag \\
    &\ \ - h \sum_{k \in S} \partial_k f_S(0) \int_{[-1, 1]^{|S|}\setminus[0, 1]^{|S|}} u_k \prod_{k \in S} K(u_k) du_k
\end{align}

The third term in the final sum tends to $0$ with $h$. The same is true for the first term as
\begin{align}
    \Big| \int_{[0, 1]^{|S|}} &\left[f_S(hu) - f_S(0) - h \sum_{k \in S} u_k \partial_k f_S(0) \right] \prod_{k \in S} K(u_k) du_k \Big| \notag\\
    &\leq \int_{[0, 1]^{|S|}} \left|f_S(hu) - f_S(0) - h \sum_{k \in S} u_k \partial_k f_S(0) \right| \prod_{k \in S} K(u_k) du_k \\
    &\leq L h^2 \sum_{k \in S} \int_{[0, 1]^{|S|}} u_K^2 \prod_{k \in S} K(u_k) du_k \\
    &\leq L \sigma_K^2 h^2.
\end{align}

Now, as $\int_{[-1, 1]^{|S|}\setminus[0, 1]^{|S|}} \prod_{k \in S} K(u_k) du_k = \frac{2^{|S|} - 1}{2^{|S|}}$, we conclude that
\begin{equation}
    \lim_{h \rightarrow 0} \esperance{\hat{f}^{\text{KDE}}_{h, S}}(0) = f_S(0) \frac{2^{|S|} - 1}{2^{|S|}} \geq e^{- A|S|} \frac{2^{|S|} - 1}{2^{|S|}} > 0.
\end{equation}
\end{proof}

\paragraph{Mirror-Image KDE} To correct the boundary bias previously introduced, a solution is to add a correction to the estimator $\hat{f}_{h ,S}^{\text{KDE}}$ near the boundary of the domain of definition. Let us define three mirroring operations for a number $x \in [0, 1]$

\begin{equation}
    M^{-1}(x) = -x; \ M^{0}(x) = x; \ M^{1}(x) = 2 - x. 
\end{equation}

\begin{figure}
    \centering
    \includegraphics[scale=0.5]{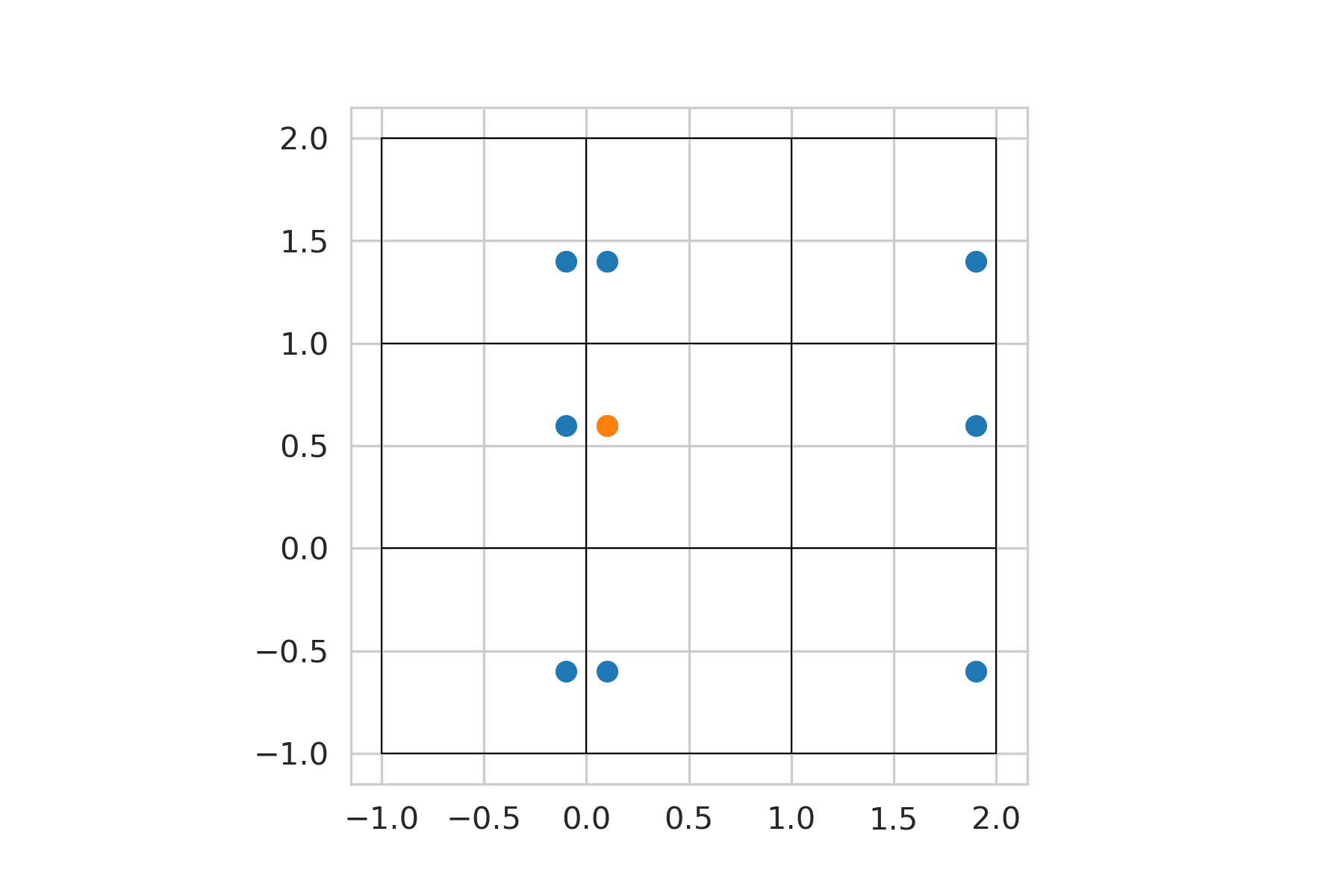}
    \caption{A datapoint in $[0, 1]^d$ (in orange) and his 8 mirror-images (in blue)}
    \label{figmirror}
\end{figure}

We define the mirror-image KDE as a KDE constructed over the sample $W$ augmented with mirror reflections of each point over all axis. An illustration of this operation in dimension $2$ is given by \cref{figmirror}.

This estimator is an extension to every dimension of the one proposed in \cite{liu2012exponential}. The formal definition is

\begin{equation}
    \hat{f}_{m ,S}^{\text{MI}}(x) = \mathbb{1}_{[0, 1]^{|S|}}(x) \frac{1}{mh^{|S|}} \sum_{i=1}^m \sum_{a \in \{-1, 0, 1\}^{|S|}} \prod_{k \in S} K\left( \frac{M^{a_k}(W_i)_k - x_k}{h} \right).
\end{equation}

\begin{figure}%
    \centering
    \subfloat[\centering A kernel is fitted over all points of the augmented dataset]{{\includegraphics[width=7cm]{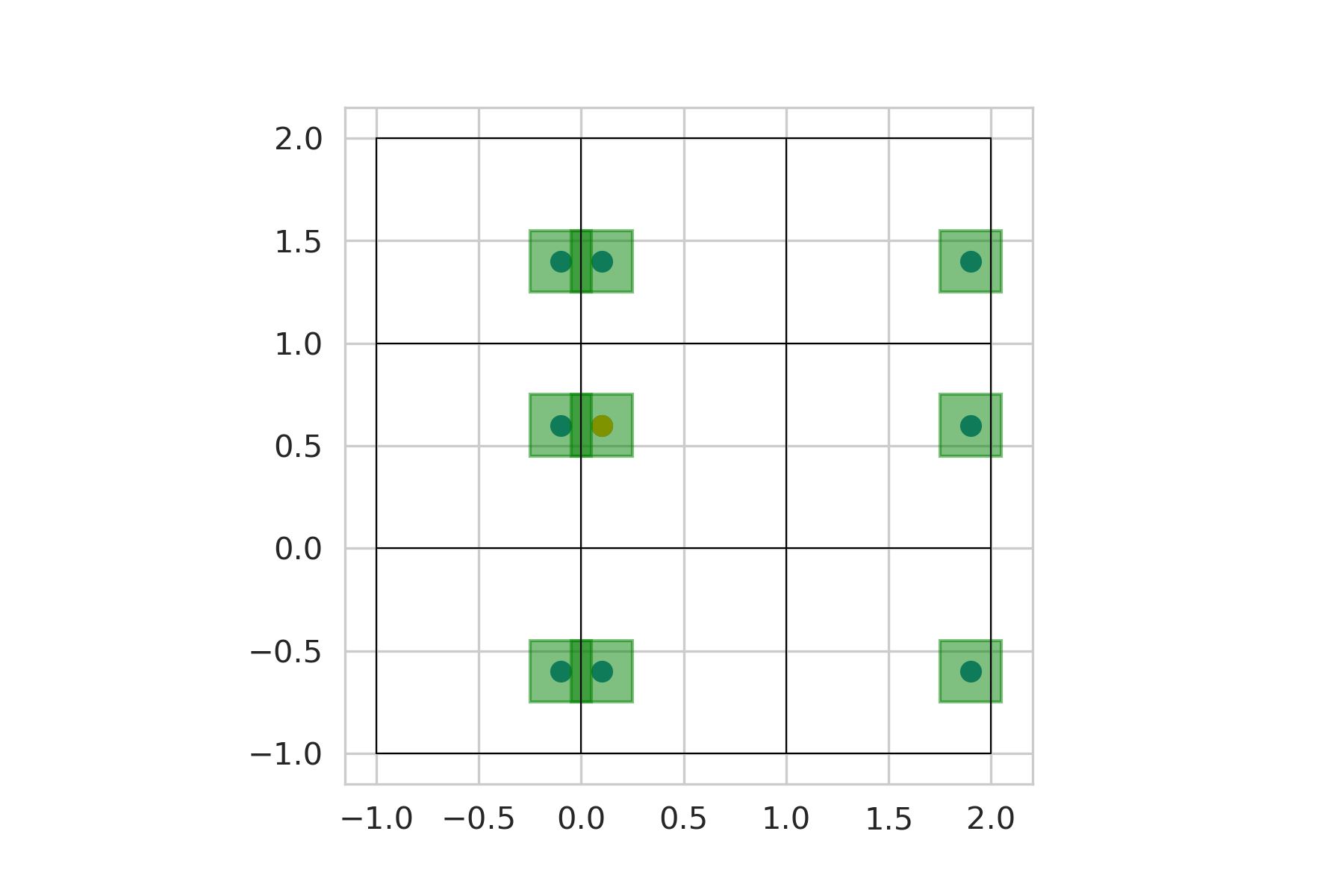} }}%
    \quad
    \subfloat[\centering Restriction to the unit hypercube]{{\includegraphics[width=7cm]{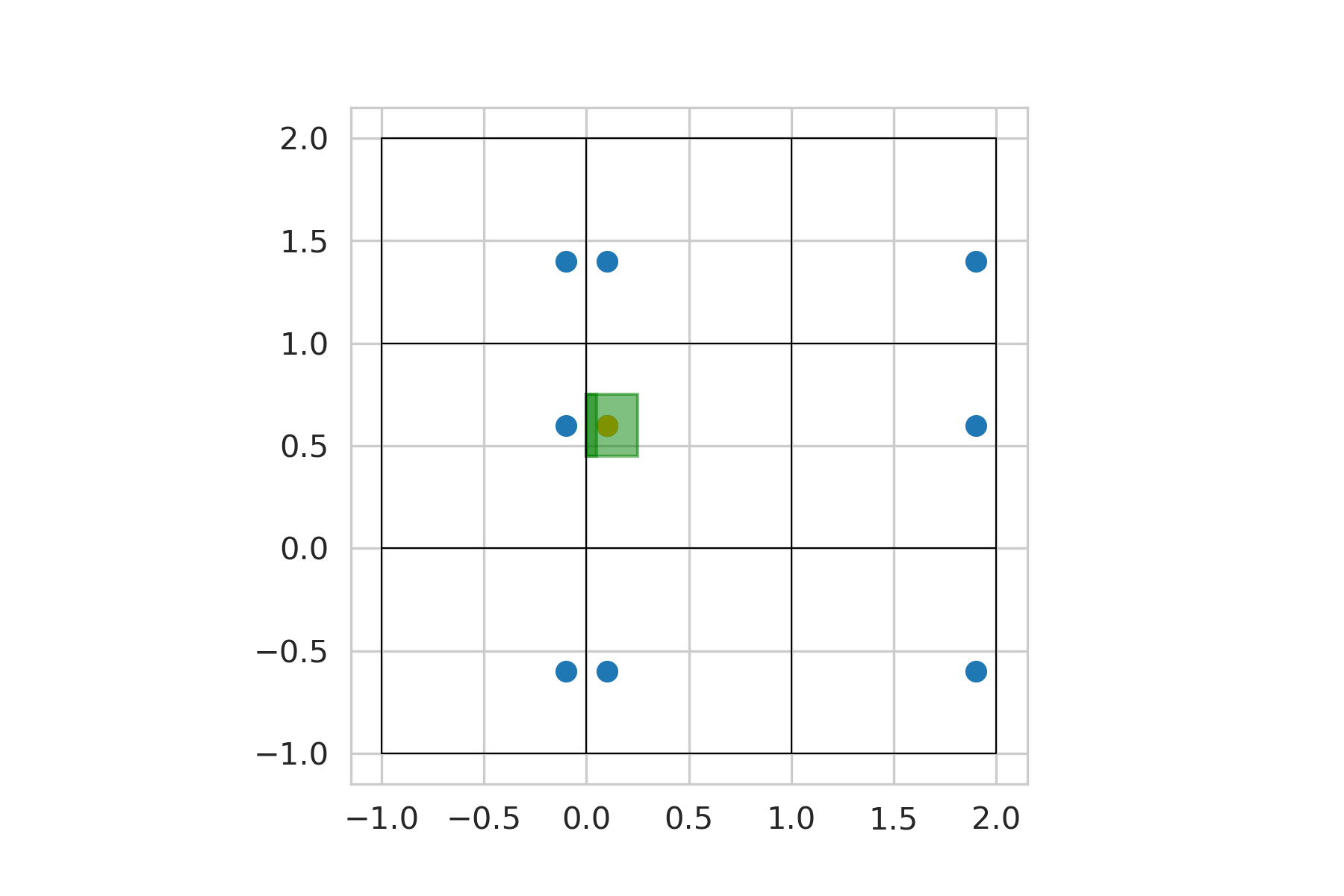} }}%
    \caption{Construction of the mirror-image KDE}%
    \label{figmiprinciple}%
\end{figure}

It consists in summing multidimensional kernel over the points of the augmented samples and restrict the domain of the obtained function to $[0, 1]^{|S|}$ as illustrated in \cref{figmiprinciple}. Roughly speaking it consists in flipping the part of $\hat{f}_{h, S}^{\text{KDE}}$ that fall outside $[0, 1]^{|S|}$ inside it. $\hat{f}_{m ,S}^{\text{MI}}$ is supported on $[0, 1]^{|S|}$ and $\int_{[0, 1]^{|S|}} \hat{f}_{m ,S}^{\text{MI}}(x) dx = 1$.

\paragraph{Bias for mirror-image KDE} Under an ad-hoc condition on the partial derivatives of $f_S$ at the boundary of $[0, 1]^{|S|}$ it is possible to bound the bias for the mirror-image KDE. Our result is an extension of the lemma 3.1 in \cite{liu2012exponential} to every dimension and every $\beta \in (0, 2]$ while the analysis in the original paper was restricted to bi-dimensional densities and $\beta = 2$. With our proof strategy, we find a better constant in the upper-bound for $|S| = 2$ and $\beta = 2$.

\begin{myprop}{Bias for mirror-image KDE}{biasmikde}
Let us assume that for all sequence $(x_n)_{n \in \mathbb{N}}$ in $[0, 1]^{|S|}$, if $x_n$ converges to a boundary point of $[0, 1]^{|S|}$, then for all $k \in S$ $\lim_{n \rightarrow \infty} \partial_kf_S(x_n) = 0$. Then
\begin{equation}
        \left\Vert f_S - \esperance{\hat{f}^{\text{MI}}_{m, S}} \right\Vert_\infty \leq C_1 h^\beta
\end{equation}
where $C_1 = L |S|^{\beta / 2} \left( 2 \Vert K \Vert_\infty \right)^{|S|}$ if $\beta < 2$ and $C_1 = L|S|$ if $\beta = 2$
\end{myprop}

\begin{proof}

    We define $f^{\text{MI}}_S$ as the function defined over $[-1, 2]^{|S|}$ such that for all $x \in [0, 1]^{|S|}$ and $a \in \{-1, 0, 1\}^{|S|}$
    \begin{equation}
        f^{\text{MI}}_S(M^{a}(x)) = f_S(X)
    \end{equation}
    where $M^{a}(x) = \left( M^{a_k}(x_k) \right)_{k \in S}$. The property that the partial derivatives of $f_S$ vanish near the boundary of $[0, 1]^{|S|}$ ensures that $\partial_k f^\text{MI}_S$ is continuous on $(-1, 2)^{|S|}$ and so $f^\text{MI}_S \in \mathcal{H}(2, L)$.
    
    Let $x \in [0, 1]^{|S|}$, we want to bound $\left| f_S(x) - \esperance{\hat{f}^{\text{MI}}_{m, S}(x)} \right|$. Assume first that $x \in [0, 1/2]^{|S|}$ and denote by $\mathcal{A}$ the set $\{ k \in S : x_k < h \}$. We start by considering the situation where $|\mathcal{A}| \geq 1$. For all $k \in \mathcal{A}$ and all $t \in [0, 1]$, $K\left( \frac{t - (2 - x_k)}{h} \right) = 0$ because the support of $K$ is $[-1, 1]$, $h \leq 1/2$ and $x_k < h$. For all $k \in S \setminus \mathcal{A}$ and all $t \in [0, 1]$, $K\left( \frac{t - (2 - x_k)}{h} \right) = 0$ and $K\left( \frac{t - (- x_k)}{h} \right) = 0$. Then the expected value of $\hat{f}^{\text{MI}}_{m, S}$ at the point $x$ can be written as
    \begin{equation}
        \esperance{\hat{f}^{\text{MI}}_{m, S}}(x) =  \sum_{\mathcal{B \subset A}} \frac{1}{h^{|S|}} \int_{[0, 1]^{|S|}} \prod_{k \in \mathcal{B}} K\left( \frac{t_k + x_k}{h} \right) \prod_{k \in S \setminus \mathcal{B}} K\left( \frac{t_k - x_k}{h} \right) f_S(t) dt
    \end{equation}
    
    Now, for $\mathcal{B} \in \mathcal{A}$, we denote $x_\mathcal{B}$ the vector such that $(x_\mathcal{B})_k = x_k$ if $k \notin \mathcal{B}$ and $(x_\mathcal{B})_k = - x_k$ if $k \in \mathcal{B}$. We have
    
    \begin{equation}
        \esperance{\hat{f}^{\text{MI}}_{m, S}}(x) =  \sum_{\mathcal{B \subset A}} \int_{\chi^S_\mathcal{B}} \prod_{k \in S} K\left( u_k \right) f_S(x_\mathcal{B} + uh) du
    \end{equation}
    
    where $\chi^S_\mathcal{B} = \left\{ u \in [-1, 1)^{|S|} : x_\mathcal{B} + uh \in [0, 1] \right\}$. We see that $\chi^S_\mathcal{B} = \prod_{k \in S} [\underline{u}_k, \bar{u}_k)$ where $\underline{u}_k = - x_k / h$ if $k \in \mathcal{B}$, $-1$ otherwise and $\bar{u}_k = - x_k / h$ if $k \in \mathcal{A} \setminus \mathcal{B}$, $1$ otherwise. What is more, as $f^{\text{MI}}_S = f_S$ on $[0, 1]^{|S|}$ we have
    \begin{equation}
        \esperance{\hat{f}^{\text{MI}}_{m, S}}(x) =  \sum_{\mathcal{B \subset A}} \int_{\chi^S_\mathcal{B}} \prod_{k \in S} K\left( u_k \right) f^{\text{MI}}_S(x_\mathcal{B} + uh) du.
    \end{equation}

    Now, as $(\chi^S_\mathcal{B})_{\mathcal{B} \subset \mathcal{A}}$ forms a partition of $[-1, 1)^{|S|}$, we have
    \begin{align}
        f_S(x) &= \sum_{\mathcal{B} \subset \mathcal{A}} \int_{\chi^S_\mathcal{B}}  \prod_{k \in S} K\left( u_k \right) f_S(x) du \\
        &= \sum_{\mathcal{B} \subset \mathcal{A}} \int_{\chi^S_\mathcal{B}}  \prod_{k \in S} K\left( u_k \right) f^\text{MI}_S(x_\mathcal{B}) du
    \end{align}
    
    We denote by $\delta_\mathcal{B}(u, \beta)$ the quantity
    
    \begin{equation}
    \left\{
        \begin{array}{ll}
            f^{\text{MI}}_S(x_\mathcal{B} + uh) - f^{\text{MI}}_S(x_\mathcal{B}) & \text{if $\beta \in (0, 1]$}  \\
            f^{\text{MI}}_S(x_\mathcal{B} + uh) - f^{\text{MI}}_S(x_\mathcal{B}) - h \sum_{k \in S} u_k \partial_k f^{\text{MI}}_S(x_\mathcal{B})  & \text{if $\beta \in (1, 2]$}
        \end{array}
        \right.
    \end{equation}
    
    From lemma \ref{lem:controlforholder2} we have
    \begin{equation}
        | \delta_\mathcal{B}(u, \beta) | \leq L h^\beta \Vert u \Vert^\beta
    \end{equation}
    
    And, as $\int xK(x) = 0$, we have
    \begin{align}
         \left| f_S(x) - \esperance{\hat{f}^{\text{MI}}_{m, S}}(x) \right| &= \left| \sum_{\mathcal{B} \subset \mathcal{A}} \int_{\chi^S_\mathcal{B}} \prod_{k \in S} K\left( u_k \right) \delta_\mathcal{B}(u) du \right| \\
         &\leq Lh^\beta \int_{[-1, 1]^{|S|}} \prod_{k \in S} K(u_k) \Vert u \Vert^\beta du
    \end{align}
    
    If $\beta = 2$
    \begin{align}
        \int_{[-1, 1]^{|S|}} \prod_{k \in S} K(u_k) \Vert u \Vert^\beta du &= \int_{[-1, 1]^{|S|}} \prod_{k \in S} K(u_k) \sum_{k \in S} u_k^2 du \\
        &= \sum_{k \in S} \int_{-1}^1 K(u) u^2 du \\
        &\leq \sum_{k \in S} \int_{-1}^1 K(u) du = |S|.
    \end{align}
    
    If $\beta < 2$
    \begin{align}
        \int_{[-1, 1]^{|S|}} \prod_{k \in S} K(u_k) \Vert u \Vert^\beta du &\leq \Vert K \Vert_\infty^{|S|} \int_{[-1, 1]^{|S|}}  \Vert u \Vert^\beta du \\
        &\leq \Vert_\infty^{|S|} \sqrt{|S|}^\beta \int_{[-1, 1]^{|S|}}  du  \\
        &= |S|^{\beta / 2} \left( 2 \Vert K \Vert_\infty \right)^{|S|}.
    \end{align}
    
    Then $\sup_{x \in [0, 1/2]^{|S|}} \left| f_S(x) - \esperance{\hat{f}^{\text{MI}}_{m, S}}(x) \right| \leq C_1 h^\beta$. By symmetry the same inequality is true when the sup is taken over $[0, 1]^{|S|}$.

\end{proof}

Then considering the mirror-image KDE leads to a correction of the boundary issue previously mentioned.

\subsubsection{Variance upper-bound}

To upper-bound the variance of the mirror-image KDE, we will use corollary 15 of \cite{kim2019uniform}. Our setting is not the same as in this paper as we deal with mirror-image KDE. Then in order to obtain the same result, we must ensure that the family of functions

\begin{equation}
    \mathcal{F}_{S}^\text{MI} = \left\{ K_{x, h}^{\text{MI}}  | x \in [0, 1]^{|S|}, h \in (0, 1/2) \right\}
\end{equation}

where

\begin{equation}
    K_{x, h}^{\text{MI}} : u \mapsto \frac{1}{h^{|S|}} \mathbb{1}_{[0, 1]^{|S|}} \sum_{a \in \{-1, 0, 1\}^{|S|}} \prod_{k \in S} K\left( \frac{M^{a_k}(u_k) - x_k}{h} \right)
\end{equation}

is a bounded VC class of function. We know that $\mathcal{F}_S$ is a bounded VC class of function. The results of section 2.6 of \cite{van1996weak} indicate that a family of functions is a bounded VC class if and only if the associated collection of sublevels is a VC class of sets. Now, we remark that the sublevels of functions in $\mathcal{F}_{S}^\text{MI}$ can be written as intersections of sublevels of functions in $\mathcal{F}_S$ intersected with $[0, 1]^{|S|}$. Then, as intersections preserve the VC class property for collection of sets (see \cite{van2009note}), $\mathcal{F}_{S}^\text{MI}$ is a bounded VC class of functions, and the corollary 15 of \cite{kim2019uniform} applies, leading to the following result.

\begin{myprop}{Variance}{variancenormeinf}

    Let $h_{m, S}$ be a bandwidth in $(0, 1/2)$ and $\delta_m \in (0, 1)$.  With probability $1 - \delta_m$
    \begin{equation}
        \left\Vert \hat{f}_{S, {h_{m, S}}} - \esperance{\hat{f}_{m, {h_{m, S}}}} \right\Vert_\infty \leq C_2 \sqrt{\frac{\log\left( 1 / h_{m, S}\right) + \log(2 / \delta_m)}{m h_{m, S}^{|S|}}}.
    \end{equation}
    The constant $C_2$ depends on $|S|$, on $\Vert K \Vert_\infty$ and on $\Vert K' \Vert_\infty$.
\end{myprop}
 




\subsubsection{Conclusion}

Now, as we have for a control of the bias and the variance term for every bandwidth $h_{m, S} \in (0, 1/2)$, by choosing appropriately $h_{m, S}$ it is possible to bound $\Vert f_S -  \hat{f}_{S, {h_{m, S}}} \Vert_\infty$.

\begin{myprop}{Convergence}{uniformconvergece}

Choosing $h_{m, S} \asymp \left(1 / m \right)^\frac{1}{2 \beta + |S|}$, it exists a constant $C_S$ such that with probability at least $1 - \delta_m$
\begin{equation}
    \Vert f_S -  \hat{f}_{S, {h_{m, S}}} \Vert_\infty \leq C_S \sqrt{\log m + 2\log\left( 2 / \delta_m \right)} \left( \frac{1}{m} \right)^{\frac{\beta}{2 \beta + |S|}}.
\end{equation}
\end{myprop}

\subsection{Uniformity over $\Sdk$}

We have just established a control in high probability for the quantity $\Vert f_S -  \hat{f}_{S, {h_{m, S}}} \Vert_\infty$ for a given $S$. Our objective is to have such a control uniformly over $\Sdk$. Applying a union-bound, we obtain the following result.

\begin{myprop}{Uniform control of uniform error over all subsets}{uniformconvergenceoversets}
Let us denote $C_k = \max_{S \in \Sdk} C_S$. We have, with probability at least $1 - S_d^k \delta_m$
\begin{equation}
    \sup_{S \in \Sdk} \left\Vert \hat{f}_{S, h^S_m} - f_S \right\Vert_\infty \leq C_k \sqrt{\log m + 2\log\left( 2 / \delta_m \right)} \left( \frac{1}{m} \right)^{\frac{2}{4 + k}}.
\end{equation}

In particular with the choice $\delta_m = \frac{2}{m S_d^k}$, we have that with probability $1 - \frac{2}{m}$
\begin{equation}
    \sup_{S \in \Sdk} \left\Vert \hat{f}_{S, h^S_m} - f_S \right\Vert_\infty \leq C_k \sqrt{3 \log m + 2\log S_d^k} \left( \frac{1}{m} \right)^{\frac{2}{4 + k}}.
\end{equation}

\end{myprop}

\section{THEOREM}\label{sec:theorem}

Let now $m_0$ be the smallest integer $m$ such that 
\begin{equation}
    C_k \sqrt{\log m + 2\log\left( 2 / \delta_m \right)} \left( \frac{1}{m} \right)^{\frac{2}{4 + k}} \leq e^{-A|S|}\left(1 - e^{-A|S|}\right).
\end{equation}

If $m \geq m_0$, we know that on an event $\mathcal{A}_m^k$ of probability at least $1 - S_d^k\delta_m$
\begin{equation}
    \Vert f_S - \hat{f}_{S, h_m} \Vert_\infty \leq C_k \sqrt{\log m + 2\log\left( 2 / \delta_m \right)} \left( \frac{1}{m} \right)^{\frac{\beta}{2\beta + k}}.
\end{equation}

Then, on $\mathcal{A}_m^k$  \eqref{unifconvergence} is satisfied with $\epsilon_S = C_k \sqrt{\log m + 2\log\left( 2 / \delta_m \right)} \left( \frac{1}{m} \right)^{\frac{\beta}{2 \beta + k}}$ for all $S \in \Sdk$. As a consequence, using proposition \ref{prop:csqunifconvergence}, we have
\begin{equation}
    \sum_{S \in \Pstar} \KL{f_S}{\hat{f}_{S, h_m^s}} \leq e^{2Ak} |\Pstar| C_k \sqrt{\log m + 2\log\left( 2 / \delta_m \right)} \left( \frac{1}{m} \right)^{\frac{2}{4 + k}}.
\end{equation}
And, on $\mathcal{A}_m^k$, for $\delta_n \in (0, 1/S_d^k)$ with probability $1 - S_d^k\delta_n$
\begin{equation}
    (P - P_n)(\log \hat{f}_\Ptilde - \log \hat{f}_\Phat) \leq 2d \sqrt{\frac{Ak}{n}} \sqrt{\log(2 / \delta_n)}.
\end{equation}

Now, with the choice $\delta_m = 1 / (2 S_d^k m)$ and $\delta_n = 1 / (2 S_d^k n)$ we obtain the following result

\begin{mythm}{Final bound}{final}
    If for all $S \in \Sdk$, $f \in \mathcal{H}(2, L)$, and for all sequence $(x_n)_{n \in \mathbb{N}}$ such that $\lim_{n \rightarrow \infty} x_n = x^*$ where $x^*$ belongs to the boundary of $[0, 1]^{|S|}$ and for all $k \in S$ $\lim_{n  \rightarrow \infty} \partial_k f_S(x_n) = 0$. With the choice $\hat{f}_S = \hat{f}^{\text{MI}}_{h_{m, S}, S}$ where $h_{m, S} \approx \left( \frac{1}{m} \right)^{\frac{1}{2 + |S|}}$, we have with probability at least $(1-1/m)(1-1/n)$
    
    \begin{align}
        \KL{f}{\hat{f}_\Phat} - \KL{f}{f_\Pstar} \leq &e^{2Ak} \sqrt{2} |\Pstar| C_k \sqrt{ \log m + \log\left( S_d^k \right)} \left( \frac{1}{m} \right)^{\frac{\beta}{2 \beta + k}} \notag\\
        &+ 2d \sqrt{\log n + \log\left( S_d^k \right)} \sqrt{\frac{Ak}{n}}
    \end{align}
\end{mythm}

Ignoring logarithmic factors, the rate of convergence of the approximation term is $\left( \frac{1}{m} \right)^{\frac{\beta}{2\beta + k}}$. The dependence of this quantity in $k$ illustrates that ISDE tackles the curse of dimensionality for the density estimation problem under KL loss in the same spirit that \cite{rebelles2015mathbb} showed that his estimator does for the squared $L_2$ loss.

Ignoring logarithmic factors again, the rate of convergence of the selection term is $\frac{1}{\sqrt{n}}$. This is a classical rate of convergence for hold-out procedures with bounded loss (see corollary 8.8 in \cite{massart2007concentration}).

The term $\log(S_d^k)$ in the upper-bound illustrates the combinatorial complexity reduction operated by ISDE. The presence of the $\log$ of the number of hypotheses is classical for hold-out procedures with bounded loss (see again corollary 8.8 in \cite{massart2007concentration}). In our context, we have reduced the combinatorial complexity from the number of partitions $P_d^k$ to the number of subsets $S_d^k$.

\section{CONCLUSION}\label{sec:conclusion}

In this paper, we have studied the convergence properties of ISDE. In particular, we have shown that under suitable assumptions on the true density and for the mirror-image KDE as marginal density estimator, we can provide an upper-bound valid with high-probability of the quantity
\begin{equation}
    \KL{f}{\hat{f}_\Phat} - \KL{f}{f_\Pstar}.
\end{equation}

This bound highlights how ISDE tackles the curse of dimensionality and reduces the combinatorial complexity of the density estimation problem under IS compared to a brute-force approach. These results offer a theoretical validation of the empirical observations presented in \cite{pujol2022isde}.

To complete the study, it let to study how the bias term $\KL{f}{f_\Pstar}$ behaves. It is hard to give a precise statement on this quantity in a general setting. One simple situation is when $f \in \mathcal{D}_k^d$. In this case $\KL{f}{f_\Pstar} = 0$. This bias can also be explicitly evaluated in some multivariate Gaussian frameworks, see appendix \ref{sec:appendixgaussian}.

\paragraph{Acknowledgement} This work was supported by the program \href{https://www.iledefrance.fr/paris-region-phd-2021}{Paris Region Ph.D.} of \href{https://www.dim-mathinnov.fr/}{DIM Mathinnov} and was partly supported by the French ANR Chair in Artificial Intelligence TopAI - ANR-19-CHIA-0001. The author is thankful to Marc Glisse and Pascal Massart for their constructive remarks on this work.

\bibliographystyle{plain}
\bibliography{bib}

\renewcommand{\thesection}{\Alph{section}}
\setcounter{section}{0}


\include{biasgaussiancase}

\end{document}

%% file: commands.tex
\usepackage{amsmath}
\usepackage{amsfonts}
\usepackage{amssymb}
\usepackage{amsthm}
\usepackage{mathtools}
\usepackage{dsfont} 
\numberwithin{equation}{section}
\usepackage{graphicx}
\usepackage{subfig}
\usepackage{esint}
\usepackage{float}
\usepackage{enumerate}
\usepackage{hyperref}
\usepackage{cleveref}
\usepackage{neuralnetwork}
\usepackage{bbold}

\usepackage{tcolorbox}
\tcbuselibrary{theorems}

\newtcbtheorem[auto counter,number within=section]{mythm}{Theorem}%
{code={\edef\@currentlabelname{#2}}, colback=black!5,colframe=black!75, before skip=20pt, after skip = 20pt}{thm}

\newtcbtheorem[use counter from=mythm]{myprop}{Proposition}%
{code={\edef\@currentlabelname{#2}}, colback=black!5,colframe=black!75, before skip=20pt, after skip = 20pt}{prop}

\newtcbtheorem[use counter from=mythm]{mylem}{Lemma}%
{code={\edef\@currentlabelname{#2}}, colback=black!5,colframe=black!75, before skip=20pt, after skip = 20pt}{lem}

\newtcbtheorem[use counter from=mythm]{mydef}{Definition}%
{code={\edef\@currentlabelname{#2}}, colback=black!5,colframe=black!75, before skip=20pt, after skip = 20pt}{def}

\newcommand{\esperance}[1]{\mathbb{E}\left[#1\right]}
\newcommand{\proba}[1]{\mathbb{P}\left[ #1 \right]}

\newcommand{\KL}[2]{\mathrm{KL}\left( #1 \| #2 \right)}

\newcommand{\argmin}[1]{\underset{#1}{\arg \ \min \ }}

\newcommand{\Pstar}{{\mathcal{P}_*}}
\newcommand{\Phat}{{\hat{\mathcal{P}}}}
\newcommand{\Ptilde}{{\Tilde{\mathcal{P}}}}
\newcommand{\Pa}{{\mathcal{P}}}

\newcommand{\Pdk}{\mathrm{Part}_d^k}
\newcommand{\Sdk}{\mathrm{Set}_d^k}
\newcommand{\Pd}{\mathrm{Part}_d}
\newcommand{\Sd}{\mathrm{Set}_d}

\newcommand{\bigzero}{\mbox{\normalfont\Large\bfseries 0}}
\newcommand{\rvline}{\hspace*{-\arraycolsep}\vline\hspace*{-\arraycolsep}}

%% file: biasgaussiancase.tex
\section{Bias term in the multivariate Gaussian framework}\label{sec:appendixgaussian}

In this appendix, we study the bias $\KL{f}{f_\Pstar}$ in a multivariate Gaussian framework where exact computations are possibles.

\subsection{Model and notations}

\paragraph{Model} If $\Sigma$ denotes a $d \times d$ definite positive matrix, we denote by $f_\Sigma$ the density of a multivariate centered Gaussian random variable with covariance $\Sigma$ and by $\Sigma_{\mathcal{P}}$ the matrix defined as follows.

\begin{equation}
\Sigma_{\Pa}(i, j) = \left\{
    \begin{array}{ll}
        \Sigma(i, j) & \mbox{if $i$ and $j$ belongs to the same bloc in $\Pa$}\\
        0 & \mbox{else.}
    \end{array}
\right.
\end{equation}
If $S_1$ and $S_2$ are subsets of $\{1, \dots, d\}$, we denote by $\Sigma(S_1, S_2)$ the $|S_1| \times |S_2|$ submatrix matrix of $\Sigma$ where we keep the intersection of rows in $S_1$ and columns in $S_2$, to keep notations compact, we write $\Sigma(S)$ instead of $\Sigma(S, S)$
 
 For a multivariate Gaussian random variable with covariance $\Sigma$, $f_\Sigma \in \mathcal{D}_d^k$ is equivalent to the fact that it exists a permutation $\sigma$ of $\{1, \dots, d\}$ such that $P_\sigma \Sigma P_\sigma^{-1}$ is block-diagonal with blocks of size smaller than $k \times k$. For clarity, in what follows, we will always consider that this property is met with $\sigma = \text{id}$, meaning that we restrict ourselves to partitions in which each block is made of consecutive features. This does not imply a loss of generality.
 
 We now consider that $f = f_{\Sigma}$ and
 \begin{equation}
     \Sigma = \Sigma_\Pa + \epsilon
     \label{smallperturb}
 \end{equation}
 
where $\Sigma_\Pa$ is a block-diagonal covariance matrix corresponding to the independence structure $\Pa$ and $\epsilon$ is a ``small" (in a sense to be defined later) definite positive matrix. The question is how this perturbation influences the bias term. In order to answer it, we must control $\KL{f}{f_\Pa}$ for all $\Pa$ in $\Pdk$.

\subsection{Some useful lemmas}

\paragraph{Computation of KL losses} The first useful result is an explicit computation of $\KL{f_\Sigma}{f_{\Sigma_{\Pa}}}$ for any $\Pa$ in $\Pdk$.

\begin{mylem}{Exact computation of KL between two centered multivariate Gaussian}{biasgaussian}
For every $\Pa \in \Pdk$
\begin{equation}
    \KL{f_\Sigma}{f_{\Sigma_\Pa}} = \frac{1}{2}\left( \sum_{S \in \Pa} \log \det \Sigma(S) - \log  \det \Sigma \right)
\end{equation}
Or if $\lambda_1 \leq \dots \leq \lambda_d$ are the eigenvalues of $\Sigma$ and $\lambda_1^\Pa \leq \dots \leq \lambda_d^\Pa$ the eigenvalues of $\Sigma_\Pa$ 
\begin{equation}
    \KL{f_\Sigma}{f_{\Sigma_\Pa}} = \frac{1}{2} \sum_{i = 1}^d \log\left(\frac{\lambda_i^\Pa}{\lambda_i}\right)
\end{equation}
\end{mylem}

\begin{proof}

The density $f_\Sigma$ has the following expression.

\begin{equation}
    f_\Sigma(x) = \frac{1}{(2 \pi)^{d / 2} (\det \Sigma)^{1 / 2}} \exp\left( - \frac{1}{2} x^\mathrm{T}\Sigma^{-1}x \right).
\end{equation}

We compute the KL divergence between $f_\Sigma$ and $f_{\Sigma_\Pa}$

\begin{align}
    \KL{f_\Sigma}{f_{\Sigma_\Pa}} &= \int \log\left( \frac{f_\Sigma(x)}{f_{\Sigma_\Pa}(x)} \right) f_\Sigma(x) dx \\
    &= \frac{1}{2} \log\frac{\det \Sigma_\Pa}{\det \Sigma} \underbrace{\int f_\Sigma(x)dx}_{=1}  \notag\\
    &\ \ +\frac{1}{2} \underbrace{\int x^\mathrm{T}\Sigma_\Pa^{-1}x f_\Sigma(x)dx}_{= \mathrm{Tr}(\Sigma_\Pa^{-1}\Sigma)}  \notag\\
    &\ \ + \frac{1}{2} \underbrace{\int x^\mathrm{T}\Sigma^{-1}x f_\Sigma(x)dx}_{= \mathrm{Tr}(\Sigma^{-1}\Sigma) = d} \\
    &= \frac{1}{2}\left( \log \det \Sigma_\Pa - \log \det \Sigma + \mathrm{Tr}\left( \Sigma_\Pa^{-1} \Sigma \right) - d \right)
\end{align}

Now, we define a permutation $\sigma_\Pa$ of $\{1, \dots, d\}$ such that :
\begin{equation}
    \Sigma_\Pa = P_{\sigma_\Pa} \begin{pmatrix}
      \Sigma(S_1)
      & \rvline & \bigzero & \rvline & \dots & \rvline & \bigzero \\
    \hline
      \bigzero & \rvline &
      \Sigma(S_2)
      & \rvline & \ddots & \rvline & \vdots \\
      \hline
      \vdots & \rvline & \ddots & \rvline & \ddots & \rvline & \bigzero \\
      \hline
      \bigzero & \rvline & \dots & \rvline  & \bigzero & \rvline & \Sigma(S_M)
    \end{pmatrix} P_{\sigma_\Pa}^{-1}
\end{equation}
where $\{S_1, \dots, S_M\}$ denotes the blocks of $\Pa$. It is then clear that $\log \det \Sigma_\Pa = \sum_{S \in \Pa} \log \det \Sigma(S)$. We also have
\begin{equation}
    \Sigma = P_{\sigma_\Pa} \begin{pmatrix}
      \Sigma(S_1)
      & \rvline & \Sigma(S_1, S_2) & \rvline & \dots & \rvline & \Sigma(S_1, S_M) \\
    \hline
      \Sigma(S_2, S_1) & \rvline &
      \Sigma(S_2)
      & \rvline & \ddots & \rvline & \vdots \\
      \hline
      \vdots & \rvline & \ddots & \rvline & \ddots & \rvline & \Sigma(S_{M-1}, S_M)  \\
      \hline
      \Sigma(S_M, S_1) & \rvline & \dots & \rvline  & \Sigma(S_{M}, S_{M-1}) & \rvline & \Sigma(S_M)
    \end{pmatrix} P_{\sigma_\Pa}^{-1}.
\end{equation}
Then
\begin{equation}
    \Sigma_\Pa^{-1}\Sigma = P_{\sigma_\Pa} \begin{pmatrix}
      I_{|S_1|}
      & \rvline & X_{1, 2} & \rvline & \dots & \rvline & X_{1, M} \\
    \hline
      X_{2, 1} & \rvline &
      I_{|S_2|}
      & \rvline & \ddots & \rvline & \vdots \\
      \hline
      \vdots & \rvline & \ddots & \rvline & \ddots & \rvline & X_{M-1, M}  \\
      \hline
      X_{M, 1} & \rvline & \dots & \rvline  & \Sigma(S_{M}, S_{M-1}) & \rvline & I_{|S_M|}
    \end{pmatrix} P_{\sigma_\Pa}^{-1}
\end{equation}
where for $i \neq j$, $X_{i, j}$ is a $|S_i| \times |S_j|$ matrix. Then $\mathrm{Tr}\left( \Sigma_\Pa^{-1} \Sigma \right) = d$.

The formulation of the result involving the eigenvalues comes from the fact that $\det \Sigma = \prod_{i=1}^d \lambda_i$ and $\det \Sigma_\Pa = \prod_{i=1}^d \lambda_i^\Pa$.
\end{proof}

\paragraph{Some computation of determinants} We define the $p\times p$ matrix
\begin{equation}
A^p_\sigma = \begin{pmatrix}
1 & \sigma & \dots & \sigma \\
\sigma & \ddots & \ddots & \vdots \\
\vdots & \ddots & \ddots & \sigma \\
\sigma & \dots & \sigma & 1
\end{pmatrix}.
\end{equation}

If $k$ divides $d$ we define

\begin{equation}
    \Sigma^{(d, k)}_\sigma = \left(\begin{array}{ c | c | c | c }
    A^k_\sigma & \mathbf{0} & \dots & \mathbf{0} \\
    \hline
    \mathbf{0} & \ddots & \ddots & \vdots \\
    \hline
     \vdots & \ddots &  \ddots & \mathbf{0} \\
     \hline
     \mathbf{0} & \dots & \mathbf{0} & A^k_\sigma
  \end{array}\right)
\end{equation}

For $\epsilon > 0$ we define
\begin{equation}
    \Sigma^{(d, k)}_{\sigma, \epsilon} = \left(\begin{array}{ c | c | c | c }
    A^k_\sigma & \mathbf{\epsilon} & \dots & \mathbf{\epsilon} \\
    \hline
    \mathbf{\epsilon} & \ddots & \ddots & \vdots \\
    \hline
     \vdots & \ddots &  \ddots & \mathbf{\epsilon} \\
     \hline
     \mathbf{\epsilon} & \dots & \mathbf{\epsilon} & A^k_\sigma
  \end{array}\right)
\end{equation}

\begin{mylem}{Property for block matrices}{propblockmat}
\begin{enumerate}[i]
    \item $\det A^p_\sigma = [1 - \sigma]^{p - 1} \left[1 + (p - 1)\sigma \right]$
    \item $\det \Sigma^{(d, k)}_{\sigma} = [1 - \sigma]^{\frac{d}{k}(k-1)}[1 + (k-1)\sigma]^{\frac{d}{k}}$
    \item $\det \Sigma^{(d, k)}_{\sigma, \epsilon} = [1 - \sigma]^{\frac{d}{k}(k-1)}[1 + (k-1)\sigma + (d - k) \epsilon] [1 + (k-1)\sigma - k\epsilon]^{\frac{d}{k} - 1}$
\end{enumerate}
\end{mylem}

\begin{proof}
(i) We start by computing the eigenvalues of $A^p_\sigma$. We remark that
\begin{equation}
    A - (1 - \sigma) I = \begin{pmatrix} \sigma & \dots & \sigma \\
                                         \vdots & & \vdots \\
                                         \sigma & \dots & \sigma \end{pmatrix}.
\end{equation}

Then it is clear that $x \in \mathbb{R}^p \in \ker\left( A^p_\sigma - (1 - \sigma) I \right) \Leftrightarrow x \in \{ y \in \mathbb{R}^p: \sum_{i=1}^p y_i = 0 \}$, which is a linear subspace of $\mathbb{R}^p$ of dimension $p-1$.

Then we remark that if $\mathbb{1}_p = \begin{pmatrix} 1 \\ \vdots \\ 1 \end{pmatrix}$ then
\begin{equation}
    A^p_\sigma \mathbb{1}_p = \left(1 + (p - 1)\sigma\right) \mathbb{1}_p
\end{equation}
Then $1 + (p - 1)\sigma$ is an eigenvalue of $A$ and it could not be of multiplicity greater than $1$ as we have just proven than $(1 - \sigma)$ has a multiplicity of $p-1$. Using the fact that the determinant is the product of the eigenvalues, we obtain
\begin{equation}
    \det A^p_\sigma = (1 - \sigma)^{p - 1} \left(1 + (p - 1)\sigma \right)
\end{equation}

The proof of (ii) follows immediately as the determinant of a block-diagonal matrix is the product of the derminants of the diagonal blocks.

(iii) We determine the eigenvalues of $\Sigma^{(d, k)}_{\sigma, \epsilon}$. To this end we will find a set of $d$ linearly independent eigenvectors. We remark that
\begin{equation}
    \Sigma^{(d, k)}_{\sigma, \epsilon} \mathbb{1}_d = (1 + (k-1) \sigma + (d - k) \epsilon) \mathbb{1}_d.
\end{equation}
Then $1 + (k-1) \sigma + (d - k) \epsilon$ is an eigenvalue of $\Sigma^{(d, k)}_{\sigma, \epsilon}$ with multiplicity at least $1$. Now, we remark that for all integer $0 \leq i \leq \frac{d}{k} - 1$ and $2 \leq j \leq k$ we have
\begin{equation}
    \Sigma^{(d, k)}_{\sigma, \epsilon} (e_{ik+1} - e_{ik+j}) = (1 - \sigma) (e_{ik+1} - e_{ik+j}).
\end{equation}
Then $(1 - \sigma)$ is an eigenvalue of $\Sigma^{(d, k)}_{\sigma, \epsilon}$ with multiplicity at least $\frac{d}{k}(k - 1)$. Finally, if for $i<j$ we denote by $\mathbb{1}^i_j = \sum_{k = i}^j e_k$, we remark that for all integer $1 \leq i \leq \frac{d}{k} - 1$
\begin{equation}
    \Sigma^{(d, k)}_{\sigma, \epsilon} \left(\mathbb{1}^{0}_{k} - \mathbb{1}^{ik+1}_{(i+1)k}\right) = \left(1 + (k-1)\sigma - k \epsilon \right) \left(\mathbb{1}^{0}_{k} - \mathbb{1}^{ik+1}_{(i+1)k}\right).
\end{equation}
Then $1 + (k-1) \sigma + (d - k) \epsilon$ is an eigenvalue of $\Sigma^{(d, k)}_{\sigma, \epsilon}$ with multiplicity at least $\frac{d}{k} - 1$.

As $1 + \frac{d}{k}(k - 1) + \left( \frac{d}{k} - 1 \right) = d$, we now that the eigenvalues of $\Sigma^{(d, k)}_{\sigma, \epsilon}$ are $1 + (k-1) \sigma + (d - k) \epsilon)$, $1 - \sigma$ and $1 + (k-1) \sigma + (d - k) \epsilon$ with multiplicity $1$, $\frac{d}{k}(k - 1)$ and $\frac{d}{k} - 1$.

\end{proof}

\subsection{Control of the bias}

\paragraph{Almost independence structure} the following property precise the KL loss between $f_{\Sigma^{(d, k)}_{\sigma, \epsilon}}$ and $f_{\Sigma^{(d, k)}_{\sigma}}$, with a particular look at the situation where $\epsilon \rightarrow 0$.
\begin{myprop}{Almost independence}{almostindep}
\begin{align}
    \KL{f_{\Sigma^{(d, k)}_{\sigma, \epsilon}}}{f_{\Sigma^{(d, k)}_{\sigma}}} = &- \frac{1}{2}\log \left( 1 + \frac{d-k}{1 + (k-1) \sigma} \epsilon\right) \notag\\
    &- \frac{1}{2} \left( \frac{d}{k} - 1 \right) \log \left( 1 - \frac{k}{1 + (k-1)\sigma} \epsilon \right)
\end{align}

At the limit $\epsilon \rightarrow 0$
\begin{equation}
    \KL{f_{\Sigma^{(d, k)}_{\sigma, \epsilon}}}{f_{\Sigma^{(d, k)}_{\sigma}}} \underset{\epsilon \rightarrow 0}{=} \frac{d(d-k)}{4(1 + (k-1)\sigma)^2} \epsilon^2 + o(\epsilon^2) 
\end{equation}
\end{myprop}

\begin{proof}
As $\Sigma^{(d, k)}_{\sigma}$ is a block-diagonal submatrix of $\Sigma^{(d, k)}_{\sigma, \epsilon}$, using lemma \ref{lem:biasgaussian} we have
\begin{equation}
    \KL{\Sigma^{(d, k)}_{\sigma, \epsilon}}{\Sigma^{(d, k)}_{\sigma}} = \frac{1}{2} \log\left( \frac{\det \Sigma^{(d, k)}_{\sigma}}{\det \Sigma^{(d, k)}_{\sigma, \epsilon}} \right)
\end{equation}

Let $\beta = 1 + (k-1)\sigma$ Now, using lemma \ref{lem:propblockmat}, we have
\begin{align}
    \KL{\Sigma^{(d, k)}_{\sigma, \epsilon}}{\Sigma^{(d, k)}_{\sigma}} &= \frac{1}{2} \left[ \log\left(\frac{\beta^{\frac{d}{k}}}{[\beta + (d - k) \epsilon] [\beta - k\epsilon]^{\frac{d}{k} - 1}}\right)\right] \\
    &= \frac{1}{2}\left[ \frac{d}{k} \log \beta - \log( \beta + (d-k)\epsilon) - \left( \frac{d}{k} - 1 \right) \log \left( \beta - k \epsilon \right) \right] \\
    &= \frac{d}{2k} \log \beta - \frac{\log \beta}{2} - \frac{1}{2}\log\left(1 + \frac{d-k}{\beta}\epsilon\right)  \notag\\
    &\  - \frac{1}{2}\left( \frac{d}{k} - 1 \right) \log \beta - \frac{1}{2}\left( \frac{d}{k} - 1 \right) \log\left(1 - \frac{k}{\beta}\epsilon \right) \\
    &= - \frac{1}{2}\log\left(1 + \frac{d-k}{\beta}\epsilon\right) - \frac{1}{2}\left( \frac{d}{k} - 1 \right) \log\left(1 - \frac{k}{\beta}\epsilon \right)
\end{align}

Now, we use that
\begin{equation}
    \log(1 + \frac{d-k}{\beta} \epsilon) \underset{\epsilon \rightarrow 0}{=} \frac{d-k}{\beta} \epsilon - \frac{(d-k)^2}{2 \beta^2} \epsilon^2 + o(\epsilon^2)
\end{equation}
and 
\begin{equation}
    \log(1 - \frac{k}{\beta} \epsilon) \underset{\epsilon \rightarrow 0}{=} -\frac{k}{\beta} \epsilon - \frac{k^2}{2 \beta^2} \epsilon^2 + o(\epsilon^2).
\end{equation}
And we obtain
\begin{align}
    \KL{\Sigma^{(d, k)}_{\sigma, \epsilon}}{\Sigma^{(d, k)}_{\sigma}} &\underset{\epsilon \rightarrow 0}{=} - \frac{d - k}{2 \beta} \epsilon + \frac{(d-k)^2}{4 \beta^2} \epsilon^2 \notag\\
    &\ \ \  + \frac{\left(\frac{d}{k} - 1 \right)k}{2 \beta} \epsilon + \frac{\left( \frac{d}{k} - 1 \right) k^2}{4 \beta^2} \epsilon^2 + o(\epsilon^2) \\
    &\underset{\epsilon \rightarrow 0}{=} \frac{(d-k)^2 + kd - k^2}{4 \beta^2} \epsilon^2 + o(\epsilon^2) \\
    &\underset{\epsilon \rightarrow 0}{=} \frac{d^2 - 2kd + k^2 + kd - k^2}{4 \beta^2} \epsilon^2 + o(\epsilon^2) \\
    &\underset{\epsilon \rightarrow 0}{=} \frac{d(d-k)}{4 \beta^2} \epsilon^2 + o(\epsilon^2) 
\end{align}
\end{proof}

\paragraph{Optimal structure for small $k$} The following proposition establish that if $\Sigma = A^d_\sigma$ and $k<d$, $\Pstar$ is composed of a maximum number of blocks of size $k$.
\begin{myprop}{Optimal Structure}{optimsrtruct}

    Suppose that $\Sigma = A^d_\sigma$. A structure $s = (s_i)_{i=1}^M$ is a list of positive integer with $\sum_{i=1}^M s_i = d$. To a structure is associated a partition with blocks of consecutive features with size $s_1, \dots, s_M$. For any structure $s$ we have
    \begin{equation}
        \KL{f_\Sigma}{f_{\Sigma_s}} = \frac{1}{2} \left( \sum_{i=1}^M \log\left(\frac{1 + (s_i - 1) \sigma}{1 - \sigma} \right) - \log\left(\frac{1 + (d-1) \sigma}{1- \sigma}\right) \right).
    \end{equation}
    
    If we denote by $p$ and $r$ the only integers such that $d = pk + r$ where $r < k$, we have
    \begin{equation}
        s^* = (\underbrace{k, \dots, k}_{p \text{ times}}, r)
    \end{equation}
\end{myprop}

\begin{proof}

We combine lemma \ref{lem:biasgaussian} and lemma \ref{lem:propblockmat} to obtain
\begin{align}
    \KL{f_\Sigma}{f_{\Sigma_s}} &= \frac{1}{2} \left( \sum_{S \in \mathcal{P}} \log \det A_\sigma^{s_i} - \log \det A_\sigma^d  \right) \\
    &= \frac{1}{2} [ \sum_{S \in \mathcal{P}} (s_i - 1) \log (1 - \sigma) + \log(1 + (s_i - 1) \sigma)  \notag\\
    &\ \ \ \ \ - (d - 1) \log (1 - \sigma) -  \log(1 + (d - 1) \sigma) ] \\
    &= \frac{1}{2} \left( \underbrace{\sum_{i = 1}^M (s_i - 1)}_{= d - M} - (d - 1)  \right) \log(1 - \sigma) \notag\\
    &\ \ + \frac{1}{2} \left( \sum_{i = 1}^M \log\left( 1 + (s_i - 1)\sigma \right) - \log(1 + (d-1) \sigma)\right)\\
    &= \frac{1}{2} \left( \sum_{i = 1}^M \log\left( \frac{1 + (s_i - 1)\sigma}{1 - \sigma} \right) - \log\left(\frac{1 + (d-1) \sigma}{1 - \sigma}\right) \right)
\end{align}

Now we want to prove that the structure minimizing $\KL{f_\Sigma}{f_{\Sigma_s}}$ is $(\underbrace{k, \dots, k}_{p \text{ times}}, r)$. To do so we start by remarking that for any $s = (s_i)_{i=1}^M \neq (k, \dots, k, r)$ it exists $i \neq j$ such that $s_i \neq k$ and $s_j \neq k$. We will prove that, for our minimization problem it is always possible to find a better structure $\Tilde{s}$ with the following
\begin{enumerate}[i]
    \item if $s_i + s_j \leq k$, $\Tilde{S} = (s_k)_{k \notin \{i, j\}} \cup (s_i + s_j)$
    \item if $ \exists l > 0: s_i + s_j = k + l$, $\Tilde{S} = (s_k)_{k \notin \{i, j\}} \cup (k, l)$
\end{enumerate}

To prove (i), we start from
\begin{align}
    2\KL{f_\Sigma}{f_{\Sigma_s}} &=  \sum_{k=1}^M \log\left(\frac{1 + (s_k - 1) \sigma}{1 - \sigma} \right) - \log\left(\frac{1 + (d-1) \sigma}{1- \sigma}\right) \\
    2\KL{f_\Sigma}{f_{\Sigma_{\Tilde{s}}}} &= \sum_{k = 1, \dots,M, k \notin\{i, j\}} \log\left(\frac{1 + (s_k - 1) \sigma}{1 - \sigma} \right) + \log\left(\frac{1 + (s_i + s_j - 1) \sigma}{1 - \sigma} \right) \notag\\
    &\ \ \ \ \ - \log\left(\frac{1 + (d-1) \sigma}{1- \sigma}\right).
\end{align}
Then to prove that $\KL{f_\Sigma}{f_{\Sigma_s}} > \KL{f_\Sigma}{f_{\Sigma_{\Tilde{s}}}}$ it is sufficient to prove that for all $a, b \geq 1$, $g(a, b) > 0$ where
\begin{equation}
    g(a, b) = \log\left(\frac{1 + (a-1)\sigma}{1 - \sigma}\right) + \log\left(\frac{1 + (b-1)\sigma}{1 - \sigma}\right) - \log\left(\frac{1 + (a+b-1)\sigma}{1 - \sigma}\right).
\end{equation}

Let us start by computing $\partial_1 g(a, b)$
\begin{align}
    \partial_1 g(a, b) &= \frac{\sigma}{1 + (a-1)\sigma} - \frac{\sigma}{1 + (a+b-1)\sigma} \\
    &= \frac{\sigma}{\left(1 + (a-1)\sigma\right) \left(1 + (a+b-1)\sigma\right)}(1 + (a+b-1)\sigma - 1 - (a-1)\sigma ) \\
    &= \frac{b \sigma^2}{\left(1 + (a-1)\sigma\right) \left(1 + (a+b-1)\sigma\right)} \geq 0.
\end{align}
Then for any $b \geq 1$, $g(a, b)$ is nondecreasing in $a$. As $a$ and $b$ play similar roles in $g(a, b)$, we have that for any $a \geq 1$, $g(a, b)$ is nondecreasing in $b$. To prove that $g(a, b) \geq 0$ it is sufficient to show that $g(1, 1) > 0$.
\begin{align}
    g(1, 1) &= \log\left(\frac{1}{1 - \sigma}\right) + \log\left(\frac{1}{1 - \sigma}\right)  - \log\left(\frac{1 + \sigma}{1 - \sigma}\right)  \\
    &= - \log\left( (1 - \sigma)(1+ \sigma) \right).
\end{align}
Now, as $\sigma \in (0, 1)$, $(1-\sigma)(1+\sigma) \in (0, 1)$. Then $\log\left((1-\sigma)(1+\sigma) \right) < 0$, implying $g(1, 1) > 0$.

To prove (ii) we start from
\begin{align}
    2\KL{f_\Sigma}{f_{\Sigma_s}} &=  \sum_{k=1}^M \log\left(\frac{1 + (s_k - 1) \sigma}{1 - \sigma} \right) - \log\left(\frac{1 + (d-1) \sigma}{1- \sigma}\right) \\
    2\KL{f_\Sigma}{f_{\Sigma_{\Tilde{s}}}} &= \sum_{k = 1, \dots,M, k \notin\{i, j\}} \log\left(\frac{1 + (s_k - 1) \sigma}{1 - \sigma} \right) + \log\left(\frac{1 + (k - 1) \sigma}{1 - \sigma} \right) \notag\\
    &\ \ \ + \log\left(\frac{1 + (l - 1) \sigma}{1 - \sigma} \right) - \log\left(\frac{1 + (d-1) \sigma}{1- \sigma}\right).
\end{align}
Then to prove that $\KL{f_\Sigma}{f_{\Sigma_s}} > \KL{f_\Sigma}{f_{\Sigma_{\Tilde{s}}}}$ it is sufficient to prove that for all $x \in [l, k]$, $h(x)$ attains its minimum at $l$ or $k$ where
\begin{equation}
    h(x) = \log\left( 1 + (x-1)\sigma \right) + \log\left( 1 + ((k+l) - x-1)\sigma \right).
\end{equation}
Let us start by computing $h'(x)$
\begin{align}
    h'(x) &= \sigma \left( \frac{1}{1 + (x-1)\sigma} - \frac{1}{1 + ((k+l)-x-1)\sigma} \right) \\
    &=\frac{\sigma\left( 1 + ((k+l)-x-1)\sigma - 1 - (x-1)\sigma \right)}{(1 + (x-1)\sigma)(1 + ((k+l)-x-1)\sigma)} \\
    &= \frac{\sigma^2}{(1 + (x-1)\sigma)(1 + ((k+l)-x-1)\sigma)}((k+l) - 2x).
\end{align}
Then $h$ increases from $l$ to $(k+l) / 2$ and decreases from $(k+l) / 2$ to $k$ and $h(l) = h(k)$ the minimum of $h$ is attained on $l$ and $k$.

\end{proof}

\paragraph{Conclusion} We finish this appendix by establishing a general upper bound of $\KL{f_\Sigma}{f_\Pstar}$ where  $\Sigma = \Sigma_{\sigma, \epsilon}^{(d, k^*)}$.

\begin{mythm}{Upper-bound for the bias in a multivariate Gaussian framework}{completegaussian}
    If $\Sigma = \Sigma_{\sigma, \epsilon}^{(d, k^*)}$ and if $k < k^*$. Let $(p, r)$ be the unique couple of integer with $0 \leq r < k$ such that $k^* = pk + r$, we have
    \begin{align}
        \KL{f_\Sigma}{f_\Pstar} \leq\  &\KL{f_{\Sigma_{\sigma, \epsilon}^{(d, k^*)}}}{f_{\Sigma_{\sigma}^{(d, k^*)}}} + \frac{dp}{2k^*} \log \left( \frac{1 + (k-1)\sigma}{1 - \sigma}\right)\notag\\
        &+ \frac{d}{2k^*} \log \left( \frac{1 + (r-1)\sigma}{1 - \sigma} \right) - \frac{d}{2k^*} \log \left( \frac{1 + (k^* - 1)\sigma}{1 - \sigma} \right)
    \end{align}
    with
    \begin{equation}
    \KL{f_{\Sigma^{(d, k)}_{\sigma, \epsilon}}}{f_{\Sigma^{(d, k)}_{\sigma}}} \underset{\epsilon \rightarrow 0}{=} \frac{d(d-k)}{4(1 + (k-1)\sigma)^2} \epsilon^2 + o(\epsilon^2) 
    \end{equation}
\end{mythm}

\begin{proof}
Let us consider
\begin{itemize}
    \item the structure $\Tilde{s} = (\underbrace{k, \dots, k}_{p \text{ times}}, r)$, and $\Pa_{\Tilde{s}}$ the associated partition of $k^*$ features,
    \item the structure $s = (\underbrace{\Tilde{s}, \dots, \Tilde{s}}_{d/k^* \text{ times}})$, and $\Pa_s$ the associated partition of $d$ features,
    \item the structure $s_0 = (\underbrace{k^*, \dots, k^*}_{d/k^* \text{ times}})$ and $\Pa_{0}$ the associated partition of $d$ features.
\end{itemize}

We can upper-bound the bias term $\KL{f_\Sigma}{f_\Pstar}$ as follows
\begin{align}
    \KL{f_\Sigma}{f_\Pstar} &\leq \KL{f_\Sigma}{f_{\Pa_s}} \\
    &= \int \log\left( \frac{f_\Sigma}{f_{\Pa_s}} \right) f_\Sigma \\
    &= \int \log\left( \frac{f_\Sigma}{f_{\Pa_0}} \right) f_\Sigma + \int \log\left( \frac{f_{\Pa_0}}{f_{\Pa_s}} \right) f_\Sigma \\
    &= \KL{\Sigma_{\sigma, \epsilon}^{(d, k^*)}}{\Sigma_{\sigma}^{(d, k^*)}} + \int \log\left( \frac{f_{\Pa_0}}{f_{\Pa_s}} \right) f_\Sigma
\end{align}

The blocks of the partition $\Pa_s$ are subsets of blocks of the partition $\Pa_0$, then
\begin{align}
    \int \log\left( \frac{f_{\Pa_0}}{f_{\Pa_s}} \right) f_\Sigma &= \int \log\left( \frac{\prod_{S \in \Pa_0} f_S}{\prod_{S \in \Pa_0} (f_{\Pa_s})_S} \right) f_\Sigma \\
    &= \sum_{S \in \Pa_0} \int \log\left( \frac{f_S}{(f_{\Pa_s})_S} \right) f_S \\
    &= \sum_{S \in \Pa_0} \KL{f_{A_\sigma^d}}{f_{(A_\sigma^d)_{\Pa_{\Tilde{s}}}}} \\
    &= \frac{d}{k^*} \KL{f_{A_\sigma^d}}{f_{(A_\sigma^d)_{\Pa_{\Tilde{s}}}}}
\end{align}
Now, using proposition \ref{prop:almostindep}

\begin{align}
    \KL{f_{A_\sigma^d}}{f_{(A_\sigma^d)_{\Pa_{\Tilde{s}}}}} &= \frac{1}{2} \bigg( p \log\left( \frac{1 + (k-1)\sigma}{1 - \sigma} \right) + \log\left( \frac{1 + (r-1)\sigma}{1 - \sigma} \right) \notag\\
    &\quad \quad \quad  - \log\left( \frac{1 + (k^*-1)\sigma}{1 - \sigma} \right) \bigg)
\end{align}
And, using proposition \ref{prop:optimsrtruct}, we have that
\begin{equation}
    \KL{f_{\Sigma^{(d, k)}_{\sigma, \epsilon}}}{f_{\Sigma^{(d, k)}_{\sigma}}} \underset{\epsilon \rightarrow 0}{=} \frac{d(d-k)}{4(1 + (k-1)\sigma)^2} \epsilon^2 + o(\epsilon^2).
\end{equation}
Then we have proven the desired upper-bound.

\end{proof}

%% file: main.bbl
\begin{thebibliography}{10}

\bibitem{gine2002rates}
Evarist Gin{\'e} and Armelle Guillou.
\newblock Rates of strong uniform consistency for multivariate kernel density
  estimators.
\newblock In {\em Annales de l'Institut Henri Poincare (B) Probability and
  Statistics}, volume~38, pages 907--921. Elsevier, 2002.

\bibitem{hall1987kullback}
Peter Hall.
\newblock On kullback-leibler loss and density estimation.
\newblock {\em The Annals of Statistics}, pages 1491--1519, 1987.

\bibitem{hasminskii1990density}
Rafael Hasminskii, Ildar Ibragimov, et~al.
\newblock On density estimation in the view of kolmogorov's ideas in
  approximation theory.
\newblock {\em The Annals of Statistics}, 18(3):999--1010, 1990.

\bibitem{kim2019uniform}
Jisu Kim, Jaehyeok Shin, Alessandro Rinaldo, and Larry Wasserman.
\newblock Uniform convergence rate of the kernel density estimator adaptive to
  intrinsic volume dimension.
\newblock In {\em International Conference on Machine Learning}, pages
  3398--3407. PMLR, 2019.

\bibitem{lepski2013multivariate}
Oleg Lepski.
\newblock Multivariate density estimation under sup-norm loss: oracle approach,
  adaptation and independence structure.
\newblock {\em Annals of Statistics}, 41(2):1005--1034, 2013.

\bibitem{liu2012exponential}
Han Liu, Larry Wasserman, and John Lafferty.
\newblock Exponential concentration for mutual information estimation with
  application to forests.
\newblock {\em Advances in Neural Information Processing Systems}, 25, 2012.

\bibitem{massart2007concentration}
Pascal Massart.
\newblock Concentration inequalities and model selection.
\newblock 2007.

\bibitem{pujol2022isde}
Louis Pujol.
\newblock Isde: Independence structure density estimation.
\newblock {\em arXiv preprint arXiv:2203.09783}, 2022.

\bibitem{rebelles2015mathbb}
Gilles Rebelles et~al.
\newblock Lp adaptive estimation of an anisotropic density under independence
  hypothesis.
\newblock {\em Electronic journal of statistics}, 9(1):106--134, 2015.

\bibitem{van2009note}
Aad Van Der~Vaart and Jon~A Wellner.
\newblock A note on bounds for vc dimensions.
\newblock {\em Institute of Mathematical Statistics collections}, 5:103, 2009.

\bibitem{van1996weak}
AW~van~der Vaart, A.W. van~der Vaart, A.~van~der Vaart, and J.~Wellner.
\newblock {\em Weak Convergence and Empirical Processes: With Applications to
  Statistics}.
\newblock Springer Series in Statistics. Springer, 1996.

\end{thebibliography}
